\documentclass[12pt]{article}
\usepackage{fullpage}
\usepackage[all]{xy}

\UseComputerModernTips
\def\diagram{\displaymath\let\labelstyle\textstyle\xymatrix@1@M=3pt@C=5pc}

\def\inline{\math\xymatrix@1@M=3pt@C=2pc}

\author{J.M. Egger\thanks{Research partially supported by NSERC}}
\title{Quillen model categories without \\ equalisers or coequalisers}
\let\implies\Rightarrow
\let\inject\mu
\let\eject\varpi
\def\of{}

\def\H{{\mathcal{H}}}
\def\K{{\mathcal{K}}} 
\def\W{{\mathcal{W}}}

\def\Smooth{{\mathcal{M}}}

\def\pair#1#2{(#1,#2)}
\def\copair#1#2{[#1,#2]}
\def\Ho#1{{\mathsf{Ho}[#1]}}

\def\Kcf{\K_{cf}}
\def\cyl#1{c} 
\def\path#1{p} 
\def\lsim{\stackrel\ell\sim}
\def\rsim{\stackrel{r}\sim}

\def\ie{{\em i.e.}}

\def\twothree{two-out-of-three}
\def\fudge{arbitrary (finite)} 
\newcommand{\id}[1][]{\mathrm{id}_{#1}}
\def\arrow#1:#2->#3;{#2 \stackrel{#1}{\longrightarrow} #3}

\def\bump{\ \par}
\def\newparagraph#1#2{\newtheorem{v#1}[blah]{#2}
  \newenvironment{#1}{\begin{v#1}\bump\em}{\end{v#1}}}
\newparagraph{defn}{Definition}
\newparagraph{defs}{Definitions}
\newparagraph{thm}{Theorem}
\newparagraph{lem}{Lemma}
\newparagraph{lems}{Lemmata}
\newparagraph{cor}{Corollary}
\newparagraph{egs}{Examples}
\newparagraph{eg}{Example}
\newparagraph{qns}{Questions}
\def\pushright#1{{
   \parfillskip=0pt            
   \widowpenalty=10000         
   \displaywidowpenalty=10000  
   \finalhyphendemerits=0      
   \leavevmode                 
   \unskip                     
   \nobreak                    
   \hfil                       
   \penalty50                  
   \hskip.2em                  
   \null                       
   \hfill                      
   {#1}                        
   \par}}                      
\def\qed{{\sc q.e.d.}}
\def\maybeskip{\addvspace{\bigskipamount}}

\newenvironment{prf}{\par\maybeskip\noindent{\bf Proof}\par}{\pushright\qed\penalty-700 \bigskip}
\newenvironment{proof}[1][]{\par\maybeskip\noindent{\bf Proof} #1\par}{\pushright\qed\penalty-700 \bigskip}

\begin{document}
  \maketitle
  \begin{abstract}
    Quillen defined a {\em model category} to be a category with
    finite limits and colimits carrying a certain extra structure.
    In this paper, we show that only finite products and coproducts
    (in addition to the certain extra structure alluded to above)
    are really necessary to construct the homotopy category.
    This leads to the interesting observation that the homotopy
    category construction could feasibly be iterated.
 \end{abstract}

\section{Introduction}

This paper concerns the definition of Quillen model category and the
most basic fact which follows from said definition. 
In particular, we show that this same result can be obtained with 
slightly weaker hypotheses; \ie, that the definition of Quillen model
category can be weakened ``at no extra cost''.  

The bone of contention is the number---or, more accurately, class---of
limits and colimits which the category is required to possess. 
We recall that a category $\K$ has \fudge\ limits if and only 
if it has equalisers and \fudge\ products, \cite[p.113]{CWM}. 
Dually, $\K$ has \fudge\ colimits if and only if it has
coequalisers and \fudge\ coproducts.

The original motivation for this result lay in the author's attempt to
apply the theory of Quillen model structures to categories arising in
the study of linear logic, \cite{Egg-Chu,Egg}; such categories
frequently possess only products and coproducts.
Understanding that such a motivation may strike a topologist as
(frankly) obscure, we conclude the article by discussing some
potential applications of this result within the more traditional
demesne of topology and geometry. 

\section{Background}

The most basic fact about Quillen model categories, alluded to above, 
is this: given a Quillen model category $\K$, one can define a
category-theoretic {\em congruence} \cite[p.52]{CWM} 
$\sim$ on a certain full subcategory of $\K$, denoted $\Kcf$, such
that the resultant quotient category $\Kcf/{\sim}$ is equivalent to
the {\em category of fractions} obtained by inverting the weak 
equivalences in $\K$.  
 
We will not repeat the definition of Quillen model structure---which 
can be found, for example, in \cite[1.1.3]{Hovey}. 
But we will quickly review the definition(s) of the homotopy
relation, $\sim$. 
For the remainder of the paper, $\K$ will denote a category with
finite products and finite coproducts equipped with a Quillen model
structure. 

\begin{defs}
  Two arrows 
  \begin{inline}{
      x \ar@<1ex>[r]^\alpha \ar@<-1ex>[r]_\beta & y }
  \end{inline}
  in $\K$ are called
  \begin{enumerate}
  \item {\em left-homotopic}, denoted $\alpha\lsim\beta$, if there
    exists a diagram 
    \begin{diagram}{ 
	& x+x \ar[dl]_-\nabla \ar[d]_-\inject \ar[r]^-{\copair\alpha\beta}  
	& y \\ x & {\cyl x} \ar[l]^-\sigma \ar[ur]+DL_-\omega }      
    \end{diagram}
    with $\mu$ a cofibration and $\sigma$ a weak equivalence.  
  \item {\em right-homotopic}, denoted $\alpha\rsim\beta$, if there
    exists a diagram
    \begin{diagram}{ 
	& {\path y} \ar[d]^-\eject & y \ar[l]_-\kappa \ar[dl]^-\Delta \\
	x \ar[ur]^-\psi \ar[r]_-{\pair\beta\gamma} 
	& y \times y   }
    \end{diagram}
    with $\eject$ a fibration and $\kappa$ a weak equivalence.
  \end{enumerate}
\end{defs}

The following statements are well-known, and their proofs do not use
the existence of limits and colimits other than products and
coproducts.  
[Indeed, the first two are almost tautologous.
But the third is quite interesting, as we shall see below.]

\begin{lems}
  \begin{enumerate}
  \item The relation $\lsim$ is reflexive, symmetric and satisfies 
    left-congruity---\ie, $\alpha\lsim\beta \implies
    \theta\of\alpha\lsim\theta\of\beta$, whenever this makes sense.
  \item The relation $\rsim$ is reflexive, symmetric and satisfies
    right-congruity---\ie, $\alpha\rsim\beta \implies
    \alpha\of\theta\rsim\beta\of\theta$, whenever this makes sense.
  \item The restriction of $\lsim$ to $\Kcf$ coincides with that of
  $\rsim$. 
  \end{enumerate}
\end{lems}

What remains to show is that the common restriction of $\lsim$ and
$\rsim$ to $\Kcf$, henceforth denoted $\sim$,\footnote{In this, my
  notation is slightly non-standard.  More often, one writes
  $\psi\sim\omega$ to mean that both $\psi\lsim\omega$ and
  $\psi\rsim\omega$ hold.} 
is transitive. 

It is worth noting at this stage that the usual proof of the
transitivity of $\sim$ (as found in \cite[Lemma 4]{Quillen}), 
requires the existence of pushouts but not of pullbacks.
Of course, there is a dual proof of the same fact which requires 
pullbacks but not pushouts.  
It would seem quite odd if the existence of pushouts-or-pullbacks were
necessary as well as sufficient. 

\section{Transitivity of $\sim$}

The following lemma, although stated in somewhat more general terms,
essentially amounts to the transitivity of $\sim$. 
Its proof is, in fact, an adaptation of the usual proof of that
$\lsim$ and $\rsim$ coincide on $\Kcf$. 

\begin{lem}\label{hr-main-lem}
Let $\alpha,\beta$ and $\gamma$ be parallel arrows $\arrow:x->y$, and
suppose $\alpha\lsim\beta\rsim\gamma$.
Then $y$ fibrant implies $\alpha\lsim\gamma$;
dually, $x$ cofibrant implies $\alpha\rsim\gamma$.
\end{lem}
\begin{prf}
  Suppose that $y$ is fibrant and that the relations
  $\alpha\lsim\beta\rsim\gamma$ are witnessed as follows:
  \begin{diagram}{ 
      & x+x \ar[dl]_-\nabla \ar[d]_-\inject \ar[r]^-{\copair\alpha\beta} 
      & y & {\path y} \ar[d]^-\eject & y \ar[l]_-\kappa \ar[dl]^-\Delta \\
      x & {\cyl x} \ar[l]^-\sigma \ar[ur]+DL_-\omega  
      & x \ar[ur]^-\psi \ar[r]_-{\pair\beta\gamma} & y \times y   }
  \end{diagram}
  with $\inject$ a cofibration, $\eject$ a fibration, and
  $\sigma,\kappa$ weak equivalences. 

  Let $\eject_0,\eject_1$ denote the two components of $\eject$ so
  that $\eject=\pair{\eject_0}{\eject_1}$.  
  By a standard argument, $y$ fibrant implies that 
  $\arrow\pi_0:y\times y->y;$ and $\arrow\pi_1:y\times y->y;$ are
  fibrations; hence also $\eject_0=\pi_0\of\eject$ and
  $\eject_1=\pi_1\of\eject$.    
  Moreover, $\eject_0\of\kappa=\id[y]=\eject_1\of\kappa$, so by
  \twothree, $\eject_0,\eject_1$ are also weak equivalences.

  Now we can factor $\copair\alpha\gamma$ through $\eject_1$ as
  follows: 
  \[ \copair\alpha\gamma 
  = \copair{\eject_1\of\kappa\of\alpha}{\eject_1\of\psi} 
  = \eject_1\of\copair{\kappa\of\alpha}{\psi}  \]
  and moreover, 
  \[ \eject_0\of\copair{\kappa\of\alpha}{\psi} 
  =\copair{\eject_0\of\kappa\of\alpha}{\eject_0\of\psi} 
  =\copair\alpha\beta 
  =\omega\of\inject. \]

  Hence, we have a diagram
  \begin{diagram}@R=12pt{ 
      x+x \ar[dd]_-\inject \ar[rr]^-{\copair\alpha\gamma}  
      \ar[dr]_-{\copair{\kappa\of\alpha}\psi} && y \\
      & {\path y} \ar[dd]^-{\eject_0} \ar[ur]+DL_-{\eject_1} \\
	 {\cyl x} \ar[dr]_-\omega \\
	 & y }
  \end{diagram}
  with $\inject$ a cofibration and $\eject_0$ a trivial fibration.
  Therefore, we can find a diagonal lift 
  \begin{diagram}@R=12pt{ 
      x+x \ar[dd]_-\inject \ar[rr]^-{\copair\alpha\gamma}  
      \ar[dr]_-{\copair{\kappa\of\alpha}\psi} && y \\
      & {\path y} \ar[dd]^-{\eject_0} \ar[ur]+DL_-{\eject_1} \\
	 {\cyl x} \ar@{.>}[ur]_-\delta \ar[dr]_-\omega \\
	 & y }
  \end{diagram}
  and so the composite $\eject_1\of\delta$ witnesses
  $\alpha\lsim\gamma$. 
\end{prf}

\begin{thm}
  The relation $\sim$ is a congruence on $\Kcf$; moreover, 
  the quotient category $\Kcf/{\sim}$ is equivalent to $\K[\W^{-1}]$. 
\end{thm}
\begin{proof}
  We have already established the first statement via a series of
  lemmata; the second is proven as in \cite{Hovey}.
\end{proof}

\section{Iterated Homotopy?} 

It is well known that, for an arbitrary Quillen model category $\K$,
the homotopy category $\Ho\K$ does not have arbitrary limits and
colimits. 
It does, however, inherit discrete limits and colimits from $\K$. 
Thus the result of this article shows that, if we can find a Quillen
model structure on $\H=\Ho\K$, then there is no obstruction to forming
a {further} homotopy category, $\Ho\H=\Ho{\Ho\K}$. 

This observation could be utilised in two opposite ways: 
one might try to find Quillen model structures on already known
homotopy categories---this might prove a useful way of studying
individual homotopy invariants; or, perhaps more interestingly, one
might attempt to `subdivide' ordinary homotopy into smaller, and
hopefully more tractable, steps.  
Let us illustrate the latter idea with an example.

Consider the concept of {\em thin homotopy} which arises in
differential geometry, \cite{Thin}.  
For our present purposes, let us define a {\em smooth space} to be a
finite disjoint union of finite-dimensional differential manifolds,
and a {\em smooth map} to be a map between smooth spaces whose
restriction to each connected component of the domain is smooth.
Then the category of smooth spaces and smooth maps, $\Smooth$, has
finite products and coproducts.  

\begin{qns}
  \begin{enumerate}
  \item Does there exist a Quillen model structure on $\Smooth$
    resulting in a suitable {\em thin homotopy category},
    $\Ho\Smooth$?  
  \item If so, does there exist a Quillen model structure on 
    $\H=\Ho\Smooth$ such that $\Ho\H$ $(=\Ho{\Ho\Smooth})$ is
    equivalent to the usual homotopy category of smooth spaces?  
  \end{enumerate}
\end{qns}

Note that, while one could replace $\Smooth$ by an even larger
category---say diffeological spaces \cite{Souriau}, or Fr\"olicher
spaces \cite{CVS}---in order to avoid using the result of this
article with respect to the first question posed above, there is no 
guarantee that one could similarly avoid the result of this article 
with respect to the second. 

\bibliography{New/Bib/wqmc} 

\begin{thebibliography}{1}

\bibitem{Thin}
A.~Caetano and R.~F. Picken.
\newblock An axiomatic definition of holonomy.
\newblock {\em Internat. J. Math.}, 5(6):835--848, 1994.

\bibitem{Egg-Chu}
Jeffrey~M. Egger.
\newblock A {Q}uillen model structure for {C}hu spaces.
\newblock In {\em Proceedings of the 21st Annual Conference on Mathematical
  Foundations of Programming Semantics (MFPS XXI), 12 May 2006}, volume 155 of
  {\em Electronic Notes in Theoretical Computer Science}, pages 361--377.
  Elsevier, 2005.

\bibitem{Egg}
Jeffrey~M. Egger.
\newblock {\em Quillen model structures, $*$-autonomous categories and
  adherence spaces}.
\newblock PhD thesis, University of Ottawa, 2006.

\bibitem{CVS}
Alfred Fr{\"o}licher and Andreas Kriegl.
\newblock {\em Linear spaces and differentiation theory}.
\newblock Pure and Applied Mathematics (New York). John Wiley \& Sons Ltd.,
  Chichester, 1988.
\newblock A Wiley-Interscience Publication.

\bibitem{Hovey}
Mark Hovey.
\newblock {\em Model categories}, volume~63 of {\em Mathematical Surveys and
  Monographs}.
\newblock American Mathematical Society, Providence, RI, 1999.

\bibitem{CWM}
Saunders Mac~Lane.
\newblock {\em Categories for the working mathematician}, volume~5 of {\em
  Graduate Texts in Mathematics}.
\newblock Springer-Verlag, New York, second edition, 1998.

\bibitem{Quillen}
Daniel~G. Quillen.
\newblock {\em Homotopical algebra}.
\newblock Lecture Notes in Mathematics, No. 43. Springer-Verlag, Berlin, 1967.

\bibitem{Souriau}
J.-M. Souriau.
\newblock Groupes diff\'erentiels.
\newblock In {\em Differential geometrical methods in mathematical physics
  (Proc. Conf., Aix-en-Provence/Salamanca, 1979)}, volume 836 of {\em Lecture
  Notes in Math.}, pages 91--128. Springer, Berlin, 1980.

\end{thebibliography}

\end{document}